\documentclass[11pt]{article}
\usepackage{amssymb}
\usepackage{amsthm}
\usepackage{graphicx}
\usepackage{graphics}
\usepackage{amsthm}

\input pictex.tex
\newcommand{\esp}{\hspace{0.01cm}}

\newcommand{\clo}{\mathrm{S}^1}

\theoremstyle{definition}

\newtheorem{thm}{Theorem}[section]

\newtheorem{prop}[thm]{Proposition}

\newtheorem{lem}[thm]{Lemma}

\newtheorem{rem}[thm]{Remark}

\input epsf

\begin{document}

\date{}
\author{Andr\'es Navas}

\title{Three remarks on one dimensional bi-Lipschitz conjugacies}
\maketitle

\noindent{\bf Abstract.} In this Note we deal with bi-Lipschitz homeomorphisms
conjugating actions 
by $C^r$ circle diffeomorphisms. Using an equivariant version of the classical
Gottschalk-Hedlund 
Lemma, we prove that such a homeomorphism is necessarily a $C^r$ diffeomorphism if
these actions are 
non free, minimal, and ergodic with respect to the Lebesgue measure. However, we
exhibit a large variety 
of examples showing that this is far from being true if the actions are non minimal.
This clarifies 
slightly the content of a classical result by Ghys and Tsuboi, who proved that,
roughly, $C^1$ 
conjugacies between non free $C^r$ one-dimensional dynamical systems are 
automatically of 
class $C^r$. All the results of this Note are contained in \cite{imca}.

\vspace{0.8cm}

\noindent{\Large {\bf Introduction}}

\vspace{0.5cm}

Let $\theta_1$ and $\theta_2$ be two non free 
actions of a finitely generated group $\Gamma$ by $C^r$ circle 
diffeomorphisms, where $r \geq 1$. Suppose that there exists some bi-Lipschitz
homeomorphism 
$\phi: \mathrm{S}^1 \rightarrow \mathrm{S}^1$ conjugating $\theta_1$ and $\theta_2$, 
{\em i.e.} such that the equality \hspace{0.02cm} $\phi \circ \theta_1 (g) = 
\theta_2(g) \circ \phi$ \hspace{0.02cm} holds for every $g \in \Gamma$. 
The problem we deal with in this Note is the following: under which 
conditions on $\theta_1$ (and $\theta_2$) the map $\phi$ is automatically a $C^r$ 
diffeomorphism? This is much inspired by the classical work \cite{GT} by Ghys and
Tsuboi, 
where the same question is addressed for $C^1$ conjugacies $\phi$ assuming that $r
\geq 2$. 
In that context they proved that $\phi$ is necessarily a $C^r$ 
diffeomorphism if there is no finite orbit; if there are finite orbits, then $\phi$ is 
a $C^r$ diffeomorphism restricted to the complementary set of these orbits. See also 
proposition 4.9 in \cite{DKN} for a closely related result in the $C^{1+\alpha}$ case.

For the non free case we show in this Note that the situation is quite different 
when $\phi$ is only assumed to be bi-Lipschitz: in general, if the actions are
minimal then 
$\phi$ is still smooth, but for the non minimal case there are a lot of bi-Lipschitz 
non smooth conjugacies.

\vspace{0.35cm}

\noindent{\bf Theorem A.} {\em Let $\theta_1$ and $\theta_2$ be two minimal non free
actions of 
a finitely generated group by $C^r$ circle diffeomorphisms, where $r \geq 1$. If
$\theta_1$ and 
$\theta_2$ are conjugated by a bi-Lipschitz circle homeomorphism $\phi$ and are
ergodic with 
respect to the Lebesgue measure, then $\phi$ is a $C^r$ diffeomorphism.}

\vspace{0.35cm}

The proof of this theorem uses a version of the classical Gottschalk-Hedlund Lemma for 
group actions. Although such a version does not appear in the literature, its proof 
is an easy modification of the classical one. We decided to include it here for 
the convenience of the reader and because of its simplicity and beauty.

Concerning the hypothesis of ergodicity, it is conjectured that minimal actions of
finitely 
generated groups by $C^r$ circle diffeomorphisms are always ergodic with respect to
the 
Lebesgue measure for $r \geq 2$. For $r \!\in ]1,2[$ the situation is more
complicated: 
if the action is non free then the same should be true, but there seem to be a lot of 
free minimal non ergodic actions (compare with \cite{brasil}). Finally, for $r = 1$ 
there are minimal non ergodic actions both in the free \cite{brasil} and the non free 
cases (these last ones can be constructed using the examples given in \cite{quas}).

Let us now consider the non minimal case. Note that a conjugacy of an action to
it-self 
is a map which centralizes this action. Moreover, if $\theta_1$ and $\theta_2$ are two 
actions by $C^r$ circle diffeomorphisms which are supposed {\em a priori} to be
conjugate 
by some $C^r$ diffeomorphism $\phi_0$, and if $\phi$ is any other bi-Lipschitz
homeomorphism 
conjugating them, then the bi-Lipschitz homeomorphism $\phi_0^{-1} \phi$ centralizes 
$\theta_1$. This is why it is so important to study the centralizer problem before 
dealing with the general conjugacy problem. At this level we prove the following 
result.

\vspace{0.35cm}

\noindent{\bf Theorem B.} {\em Let $\Gamma$ be any finitely generated group of $C^2$
circle 
diffeomorphisms whose action is non minimal and for which the stabilizers of points 
are either trivial or infinite cyclic. Then there exists a bi-Lipschitz circle 
homeomorphism which is not $C^1$ and which commutes with every element of $\Gamma$. 
Moreover, such a homeomorphism can be taken to be non differentiable on every 
open interval of the circle.}

\vspace{0.35cm}

The hypothesis on stabilizers is not very strong. For instance, it is always
satisfied for 
real-analytic non minimal actions without finite orbits. (This result is due to
Hector; a 
complete proof appears in the Appendix of \cite{quasim}.) Of course, it is also
satisfied 
by many other smooth non real-analytic interesting actions. Without this hypothesis it 
is easy to see that, in some cases, bi-Lipschitz conjugacies are forced to be smooth. 

We finish with an example where the conjugacy problem cannot be reduced (in a very
strong 
sense) to a problem of centralizers. It would be interesting to know if the 
examples of the following theorem can be real-analytic.

\vspace{0.35cm}

\noindent{\bf Theorem C.} {\em There exist two finitely generated groups of
$C^{\infty}$ circle 
diffeomorphisms acting non freely and without finite orbits which are bi-Lipschitz
conjugate 
but for which there is no $C^{1}$ circle diffeomorphism conjugating them.}

\vspace{0.35cm}

In what follows we will consider only orientation preserving maps, but the results
can be easily 
extended to the non orientation preserving case (we leave this as a task to the
reader). Moreover, 
by using standard methods, the results of this Note can be generalized into the
context of 
codimension one foliations or general one-dimensional pseudo-groups.

\vspace{0.4cm}

\noindent{\bf Acknowledgments.} This work was motivated by a question asked to the
author by 
\'E. Ghys, to whom I would like to extend my gratitude. I would also thank T. Tsuboi
for 
useful comments and suggestions, as well as for his invitation to the University of 
Tokyo where this Note was mostly written.


\section{The minimal case}

\subsection{A Gottschalk-Hedlund Lemma for group actions}

\hspace{0.35cm} Let $X$ be a compact metric space and $\Gamma$ 
a finitely generated group acting on it by homeomorphisms. A 
{\em cocycle} associated to this action is a map $c: \Gamma 
\times X \rightarrow \mathbb{R}$ such that for each fixed 
$f \in \Gamma$ the map \esp $x \mapsto c(f,x)$ \esp is continuous, 
and such that for every $f,g$ in $\Gamma$ and every $x \in X$ one has
\begin{equation}
c(fg,x) = c(g,x) + c(f,g(x)).
\label{coc}
\end{equation}

\vspace{0.15cm}

\begin{lem} {\em Suppose that the $\Gamma$-action on $X$ is minimal. 
Then the following are equivalent:}

\noindent{(i) \em there exists some $x_0 \in X$ and some constant $C > 0$ 
such that $|c(f,x_0)| \leq C$ for every $f \in \Gamma$,}

\noindent{(ii) \em there exists some continuous function $\varphi: X \rightarrow 
\mathbb{R}$ such that $c(f,x) = \varphi (f(x)) - \varphi(x)$ for all 
$f \in \Gamma$ and all $x \in X$.}
\end{lem}

\vspace{-0.15cm}

\noindent{\bf Proof.} If the second condition is satisfied then 
$$|c(f,x_0)| \leq |\varphi(f(x_0))| + |\varphi(x_0)|\leq 2 \| \varphi \|_{C^0},$$
which proves the validity of condition (i).

Reciprocally, let us suppose that the first condition holds. For each $f \in 
\Gamma$ consider the homeomorphism $\hat{f}$ of the space $X \times \mathbb{R}$ 
defined by $\hat{f}(x,t) = (f(x),t+c(f,x))$. It is easy to see that the cocycle 
relation (\ref{coc}) implies that this defines a group action of $\Gamma$ on 
$X \times \mathbb{R}$, in the sense that $\hat{f} \hat{g} = \widehat{fg}$ for 
all $f,g$ in $\Gamma$. Moreover, condition (i) implies that the orbit of 
the point $(x_0,0)$ under this action is bounded; in particular, its 
closure is a (non empty) compact invariant set. Using Zorn's lemma, 
one easily deduces the existence of a minimal non empty compact 
invariant subset $M$ of $X \times \mathbb{R}$. We claim that this 
subset is the graph of a continuous function from $X$ to $\mathbb{R}$.

First of all, since the action of $\Gamma$ on $X$ is minimal, the projection 
of $M$ on $X$ is the whole space. Moreover, if $(\bar{x},t_1)$ and 
$(\bar{x},t_2)$ belong to $M$ for some $\bar{x} \in X$ and 
some $t_1 \neq t_2$, then this implies that $M \cap M_t \neq 
\emptyset$, where $t = t_2 - t_1 \neq 0$ and $M_t = \{(x,s+t): \esp 
(x,s) \in M \}$. Note that the $\Gamma$-action on $X \times \mathbb{R}$ 
commutes with the map $(x,s) \mapsto (x,s+t)$; in particular, $M_t$ is also 
invariant. But since $M$ is minimal, this implies that $M = M_t$. One then 
concludes that \esp $M = M_t = M_{2t} = \ldots$, which is impossible since 
$M$ is compact.

We have then proved that for every $x \in X$ the set $M$ contains exactly 
one point of the form $(x,t)$. Putting $\varphi(x) = t$ one obtains a 
function form $X$ to $\mathbb{R}$, which is continuous, since its 
graph (which coincides with $M$) is compact. 

Finally, since the graph of $\varphi$ is invariant by the action, 
for all $f \in \Gamma$ and all $x \in X$ the point $\hat{f}(x,\varphi(x)) 
= (f(x),\varphi(x) + c(f,x))$ must be of the form $(f(x),\varphi(f(x)))$, 
which implies that \esp $c(f,x) = \varphi(f(x)) - \varphi(x)$.

\vspace{0.15cm}

\begin{lem} {\em Let $X$ be a compact metric space and $\Gamma$ a finitely generated 
group acting on it by homeomorphisms. Suppose that the $\Gamma$-action on $X$ is 
minimal and ergodic with respect to some probability measure $\mu$, and let 
$c$ be a cocycle associated to this action. If $\varphi$ is a function in 
$L^{\infty}_{\mu}(X)$ such that for all $f \in \Gamma$ and $\mu$ 
almost every $x \in X$ one has} 
\begin{equation}
c(f,x) = \varphi (f(x)) - \varphi(x),
\label{cacha}
\end{equation}
{\em then there exists some continuous function $\tilde{\varphi}: 
X \rightarrow \mathbb{R}$ which coincides $\mu$ a.e. with $\varphi$ 
and such that for all $f \in \Gamma$ and all $x \in X$ one has}
\begin{equation}
c(f,x) = \tilde{\varphi} (f(x)) - \tilde{\varphi} (x).
\label{cacha2}
\end{equation}
\label{cle}
\end{lem}

\vspace{-0.4cm}

\noindent{\bf Proof.} Let $Y_0$ be the set of points in which (\ref{cacha}) 
does not hold for some $f \in \Gamma$. Since $\Gamma$ is finitely generated, 
$\mu(Y_0) = 0$. Let $Y_1'$ the complementary set of the essential support of 
$\varphi$, and let $Y_1 = \cup_{f \in \Gamma} f(Y_1')$. Take a point $x_0$ 
in the full measure set $X \setminus (Y_0 \cup Y_1)$. Equation (\ref{cacha}) 
then gives \esp 
$|c(f,x_0)| \leq 2 \| \varphi \|_{L^\infty}$ \esp for all $f \in \Gamma$. 
By the preceding lemma, there exists 
some continuous function $\tilde{\varphi}: X \rightarrow \mathbb{R}$ such that 
(\ref{cacha2}) holds for every $x$ and $f$. This implies that $\mu$ a.e. we have 
$$\tilde{\varphi} \circ f - \tilde{\varphi} = \varphi \circ f - \varphi,$$
and so 
$$\tilde{\varphi} - \varphi = (\tilde{\varphi} - \varphi) \circ f.$$
Since the $\Gamma$-action on $X$ is assumed to be $\mu$-ergodic, the 
difference $\tilde{\varphi} - \varphi$ has to be $\mu$ a.e. constant. 
Finally, changing $\tilde{\varphi}$ by some $\tilde{\varphi} + C$, 
we may force this constant to be equal to zero.


\subsection{Proof of Theorem A}

\hspace{0.35cm} Note that if $\phi$ is a bi-Lipschitz homeomorphism of the circle 
conjugating the actions $\theta_1$ and $\theta_2$ of our group $\Gamma$, then 
$\phi$ and $\phi^{-1}$ are almost everywhere differentiable with $L^{\infty}$ 
functions as derivatives. Therefore, the function $x \mapsto \log(\phi'(x))$ 
is also $L^{\infty}$. The relation \esp $\theta_1 (f) = \phi^{-1} 
\circ \theta_2(f) \circ \phi$ \esp gives almost everywhere 
$$\log(\theta_1(f)'(x)) = \log(\phi'(x)) - \log(\phi'(\theta_1(f)(x))) 
+ \log(\theta_2(f)'(\phi(x))).$$
Putting $\varphi = - \log(\phi')$ and $c(f,x) = \log(\theta_1(f)'(x)) - 
\log(\theta_2(f)'(\phi(x)))$ this gives, for all $f \in \Gamma$ 
and almost every $x \in \mathrm{S}^1$, 
$$c(f,x) = \varphi(\theta_1(f)(x)) - \varphi(x).$$
One easily checks the cocycle relation
$$c(fg,x) = c(g,x) + c(f,\theta_1(g)(x)).$$
Since the $\theta_1$-action is supposed to be ergodic, Lemma \ref{cle} 
gives the existence of a continuous function $\tilde{\varphi}$ which 
coincides almost everywhere with $\varphi$ and such that (\ref{cacha2}) 
holds for every $x$ and $f$. By integrating, one concludes that the 
derivative of $\phi$ is well defined everywhere and coincides 
with $\exp(-\tilde{\varphi})$. In particular, $\phi$ is of 
class $C^1$, and interviewing the roles of $\theta_1$ and 
$\theta_2$, one concludes that $\phi$ is a $C^1$ diffeomorphism. 
In order to prove that $\phi$ is a $C^r$ diffeomorphism, one can 
use the main result of \cite{GT} for $r \geq 2$, as well as 
Proposition 4.4 of \cite{DKN} for the $C^{1 + \alpha}$ case. 



\section{The non minimal case}


\subsection{Non smooth bi-Lipschitz centralizers}

\hspace{0.35cm} Before passing to the proof of Theorem B, let us explain the main
idea by giving 
a very simple and general construction (which seems to be well known to the
specialists) 
of a non smooth bi-Lipschitz homeomorphism centralizing an interval diffeomorphism 
without interior fixed points.

Let $f$ be a $C^2$ diffeomorphism of $I = [a,b]$ such that $f^n(x)$ converges to $a$
as $n$ goes to 
infinity for every $x \! \in [a,b[$. Fix any point $c \!\in ]a,b[$, and consider any
bi-Lipschitz 
homeomorphism $h$ from the interval $[f(c),c]$ to itself. Extending $h$ to $]a,b[$
in such a way 
that $f h = h f$, and then putting $h(a)=a$ and $h(b)=b$, we obtain a well defined
self-homeomorphism 
of $[a,b]$ (still denoted by $h$). We claim that this globally defined $h$ is still
bi-Lipschitz. 
More precisely, if $M$ is a bi-Lipschitz constant for $h$ on $[f(c),c]$, then $M
e^V$ is a 
bi-Lipschitz constant for $h$ on $[a,b]$, where $V$ is the total variation of the
logarithm 
of the derivative of $f$:
$$V = var (\log(f')) = 
\sup_{a \leq a_0 \leq a_1 \leq \ldots \leq a_n \leq b} \sum_{i=0}^{n-1} 
\big| \log(f'(a_{i+1})) - \log(f'(a_i)) \big| = \int_{a}^{b} \left|
\frac{f''(s)}{f'(s)} 
\right| \hspace{0.02cm} ds.$$
Indeed, let us suppose for instance that $x$ belongs to $f^n([f(c),c])$ for some $n
\geq 0$, 
and that $h$ has a well defined derivative at the point $f^{-n}(x) \in [f(c),c]$
which is less or 
equal than $M$. (Note that this is the case for almost every $x \in
[f^{n+1}(c),f^n(c)]$.)  Because 
of the relation $h=f^n h f^{-n}$ one has the inequality
\begin{equation}
h'(x) = h'(f^{-n}(x)) \cdot \frac{(f^n)'(hf^{-n}(x))}{(f^n)'(f^{-n})(x)} \leq 
M \cdot \frac{(f^n)'(hf^{-n}(x))}{(f^n)'(f^{-n})(x)}.
\label{unilita}
\end{equation}
Now, putting $y = f^{-n}(x) \in [f(c),c]$ and $z = h(y) \in [f(c),c]$, we have 
$$\left| \log \Big( \frac{(f^n)'(z)}{(f^n)'(y)} \Big) \right| = 
\left| \log \Big( \frac{\prod_{i=0}^{n-1} f'(f^i(z))}{\prod_{i=0}^{n-1} f'(f^i(y))}
\Big) \right|
\leq \sum_{i=0}^{n-1} \Big| \log(f'(f^i(z))) - \log(f'(f^i(y))) \Big| \leq V.$$
Introducing this last inequality into (\ref{unilita}) one obtains $h'(x) \leq M
e^{V}$. Since 
$x$ was a generic point, this shows that $h$ has Lipschitz constant bounded by
$Me^V$. The 
very same argument can be used to check a simlar bound for the 
Lipschitz constant of $h^{-1}$. 

For the proof of Theorem B we will try to perform 
an analoous construction. For simplicity, we will give a 
complete proof only for the first claim of the theorem, leaving to the 
reader the task of adapting our arguments to prove the second (and stronger) claim 
concerning the non differentiability on every open interval for some centralizing 
bi-Lipschitz homeomorphism.

Let us start by recalling that if $\Gamma$ is group of $C^2$ circle diffeomorphisms
(and 
more generally of circle homeomorphisms) whose action is non minimal, then there are
two 
possibilities: either $\Gamma$ preserves a minimal Cantor set (called the {\em
exceptional} 
minimal set), or $\Gamma$ has finite orbits \cite{ghys-EM}. Let us consider the
first case, 
which is dynamically more interesting. Fix any connected component $]a,b[$ of the
complementary 
of the exceptional minimal set. By a result due to Hector, the stabilizer in
$\Gamma$ of 
$I = [a,b]$ is non trivial (see Lemma 2.7 in \cite{ghys-top}), and so by the
hypothesis of 
the theorem it is infinite cyclic. Fix a generator $f$ for this stabilizer. If the
restriction 
of $f$ to $I$ is trivial we let $h$ be any bi-Lipschitz non $C^1$ homeomorphism of
$I$. If 
not, fix $[\bar{a},\bar{b}] \subset [a,b]$ such that $f^n(x) \neq x$ for every $x \in 
]\bar{a},\bar{b}[$, and $f(\bar{a}) = \bar{a}$ and $f(\bar{b}) = \bar{b}$. Changing
$f$ 
by $f^{-1}$ if necessary, we may assume that $f^n(x)$ converges to $\bar{a}$ as $n$
goes 
to infinity for every $x \in [\bar{a},\bar{b}[$. As before consider any point
$\bar{c}$ 
in $]\bar{a},\bar{b}[$, and consider any bi-Lipschitz non $C^1$ homeomorphism $h$ of 
$[f(\bar{c}),\bar{c}]$. This homeomorphism extends in a unique way to a bi-Lipschitz 
homeomorphism of $[a,b]$ commuting with the restriction of $f$ to $[\bar{a},\bar{b}]$ 
and which is the identity on $I \setminus [\bar{a},\bar{b}]$. 

By the hypothesis on stabilizers, it is easy to see that there exists a unique
extension of 
$h$ into a circle homeomorphism (still denoted by $h$) which commutes with (every
element 
of) $\Gamma$ and coincides with the identity in the complementary set of $\cup_{g
\in \Gamma} 
\hspace{0.1cm} g(]a,b[)$. We claim that this extension is still bi-Lipschitz. More
precisely, 
fixing a finite system $\mathcal{G} = \{g_1,\ldots,g_k\}$ of generators of $\Gamma$,
denoting 
by $V$ the supremum for the variation of the logarithm of the derivatives of these
generators, 
and choosing a bi-Lipschitz constant $M$ for $h$ on $[a,b]$, we claim that $h$ has 
bi-Lipschitz constant smaller or equal than $Me^{kV}$ over the whole circle. The 
proof of this claim is similar to that of the case of the interval ({\em i.e.} the
one given at the beginning of this Section). Let us choose for instance a point 
$x \in \cup_{g \in \Gamma} \hspace{0.1cm} (g(I) \setminus I)$, and let's try to 
estimate $h'(x)$. To do this, let's take a minimal $n \in \mathbb{N}$ for which 
there exists some $g = g_{i_n} \circ \ldots g_{i_1} \in \Gamma$ with each $g_{i_j}$ 
belonging to $\mathcal{G}$ and such that $g(x) \in I$. The minimality of $n$ implies 
that the intervals $I, g_{i_n}^{-1}(I),g_{i_{n-1}}^{-1}g_{i_n}^{-1}(I),
\ldots,g_{i_1}^{-1} \cdots g_{i_n}^{-1}(I)$ have disjoint interiors. Using 
the relation $h = g^{-1} h g$ one obtains, for a generic $x \in g^{-1}(I)$,
\begin{equation}
h'(x) = h'(g(x)) \cdot \frac{g'(x)}{g'(h(x))} \leq M \cdot
\frac{g'(x)}{g'(y)},
\label{jime}
\end{equation}
where $y = h(x) \in g^{-1}(I)$. Then using only the fact that the total
variation for the logarithm of the derivative of each $g_i$ is bounded by
$V$, one obtains
$$\left| \log \Big( \frac{g'(x)}{g'(y)} \Big) \right|
\leq \sum_{j = 0}^{n-1} \big| \log(g_{i_{j+1}}'(g_{i_j} \cdots g_{i_1})(x))
- \log(g_{i_{j+1}}'(g_{i_j} \cdots g_{i_1})(y)) \big| \leq 
\sum_{i=1}^k var \big( \log(g_i') \big) \leq kV.$$
Therefore, from (\ref{jime}) one concludes that $h'(x) \leq M e^{kV}$, as
desired.

Let us now consider the case of finite orbits. If $\Gamma$ is finite then consider any 
bi-Lipschitz non differentiable circle homeomorphism commuting with its (finite order) 
generator. If $\Gamma$ is infinite, then because of H\"older and Denjoy Theorems the
action 
of $\Gamma$ cannot be free. Take a non trivial element $f \in \Gamma$ having fixed
points, 
and let $I$ be some connected component of the complementary set of the union of the
finite 
orbits. Note that $f$ must fix all the points of these orbits. So, proceeding as 
in the previous case with $I$ and $f$ one can construct a bi-Lipschitz non 
differentiable circle homeomorphism centralizing $\Gamma$.


\subsection{Bi-Lipschitz conjugate actions which are non $C^1$ conjugate}
 
\hspace{0.35cm} Before entering into the proof of Theorem C, we would like to insist
on the fact 
that the constructions we propose are rather artificial, and definitively it would
be much more 
interesting to give real-analytic examples of groups sharing a similar conjugacy
property.

Let us begin by considering a very simple action on the interval
illustrating the main idea. For this, let us fix a sequence
$(\ell_n)_{n \in \mathbb{Z}}$ of positive real numbers such that
$\ell_n / \ell_{n+1}$ converges to 1 as $|n|$ goes to infinity,
such that $\ell_{2n} = \ell_{2n+1}$ for every $n \in \mathbb{Z}$, and such
that $\sum_{n \in \mathbb{Z}} \ell_n = 1$.  Then define another sequence
$(\bar{\ell}_n)_{n \in \mathbb{Z}}$ by $\bar{\ell}_{2n} = 4 \ell_n /3$
and $\bar{\ell}_{2n+1} = 2 \ell_{2n+1} / 3$. Note that 
$\sum_{n \in \mathbb{Z}} \bar{\ell}_n = 1$.

For each $n \in \mathbb{Z}$ consider a diffeomorphism $f_n$ from the interval
$$I_n = \Big[ \sum_{i < n} \ell_i, \sum_{i \leq n} \ell_i \Big]$$
to it-self without interior fixed points. Let $\phi_0$ be the 
homeomorphism of $[0,1]$ whose restriction to each $I_n$
is the affine map sending $I_n$ to
$$\bar{I}_n = \Big[ \sum_{i < n} \bar{\ell}_i, \sum_{i \leq n}
\bar{\ell}_i \Big],$$
and let $\bar{f}_n$ be the diffeomorphism of $\bar{I}_n$ 
defined by $\bar{f}_n = \phi_0 f_n \phi_0^{-1}$. It is easy to see that 
if the maps $f_n$ are well chosen (for instance, if they are infinitely 
tangent to the identity at the extreme points and their $C^r$ norm 
converge to zero exponentially fast as $|n|$ goes to infinity for every 
$r \geq 2$), then the map $f$ defined by $f(x) = f_n(x)$ for every $x \in 
I_n$ and $f(0)=0$ and $f(1)=1$, as well as $\bar{f} = \phi_0 f \phi_0^{-1}$, 
are $C^{\infty}$ diffeomorphisms of $[0,1]$ which are infinitely tangent to 
the identity at the extreme points. Moreover, it follows from the definitions 
that $\phi_0$ is a bi-Lipschitz homeomorphism conjugating them. We claim 
however that there is no $C^1$ diffeomorphism conjugating $f$ and
$\bar{f}$. Indeed, for every homeomorphism $\phi$ conjugating $f$ 
and $\bar{f}$ there exists a fixed $N \in \mathbb{N}$ such that for 
$\phi(I_n) = \bar{I}_{n+N}$ for all $n \in \mathbb{Z}$. If such $\phi$ 
was of class $C^1$ then using the continuity of $\phi'$ at $1$ one could 
conclude that, as $n \rightarrow \infty$, 
$$\frac{|\bar{I}_{n+N}|}{|I_n|} \longrightarrow \phi'(1).$$
However, the left hand side expression does not converge. Indeed, if $N$ 
is even then as $n \rightarrow \infty$ one has 
$$\frac{|\bar{I}_{2n+N}|}{|I_{2n}|} \longrightarrow \frac{4}{3} \qquad 
\mbox{and} \qquad \frac{|\bar{I}_{2n+1+N}|}{|I_{2n+1}|} \longrightarrow 
\frac{2}{3},$$
whereas if $N$ is odd then as $n \rightarrow \infty$ one has
$$\frac{|\bar{I}_{2n+N}|}{|I_{2n}|} \longrightarrow \frac{2}{3} \qquad 
\mbox{and} \qquad \frac{|\bar{I}_{2n+1+N}|}{|I_{2n+1}|} \longrightarrow 
\frac{4}{3}.$$

Now in order to obtain an example with an exceptional minimal set we will try 
to ``glue'' the preceding construction in one of the connected components of the 
complement of such a minimal set. To be more precise, let us consider the injection 
$\theta: G \rightarrow \mathrm{Diff}_+^{\infty}(\clo)$ of the Thompson group $G$
obtained 
by the method of \S III.1 of \cite{GS} by using a map satisfying the properties (I), 
(II) and (III)$_{\infty}$ therein, and having an interval of fixed points. The 
corresponding action admits an exceptional minimal set, and we can fix an interval 
$I$ contained in one of the connected components $J$ of the complement of this set 
in such a way that the restriction to $I$ of every element of $G$ fixing $J$ 
coincides with the identity map. Let $\phi_I$ be the affine map sending $[0,1]$ 
to $I$, and let $h \in \mathrm{Diff}_+^{\infty}(\clo)$ (resp. $\bar{h}$) be defined by 
$h(x) = \phi_I f \phi_I^{-1}(x)$ for $x \in I$ and $h(x) = x$ for $x \notin I$
(resp. $\bar{h}(x) 
= \phi_I \bar{f} \phi_I^{-1}(x)$ for $x \in I$ and $\bar{h}(x) = x$ for $x \notin
I$). Now consider 
the induced group $\Gamma$ which is a quotient of the free product between $G$ and
$\mathbb{Z}$. 
This group has two actions $\theta_1$ and $\theta_2$ by $C^{\infty}$ circle
diffeomorphisms, 
depending if we choose $h$ or $\bar{h}$ as the generator of $\mathbb{Z}$. These
actions are 
clearly bi-Lipschitz conjugate, but as before it is easy to see that they are non
$C^1$ 
conjugate.


\begin{small}

\vspace{0.1cm}

\noindent Andr\'es Navas\\

\noindent Univ. de Santiago de Chile, Alameda 3363, Santiago, Chile 
(andnavas@uchile.cl)\\

\end{small}


\begin{thebibliography}{Dillo 83}


\bibitem{DKN} {\sc Deroin, B., Kleptsyn, V. \& Navas, A.} Sur la dynamique 
unidimensionnelle en r\'egularit\'e interm\'ediaire. To appear in 
{\em Acta Mathematica}.\\

\bibitem{ghys-EM} {\sc Ghys, \'E.} Groups acting on the circle. {\em L'Enseign. Math.} 
{\bf 47} (2001), 329-407.\\

\bibitem{ghys-top} {\sc Ghys, \'E.} Classe d'Euler et minimal exceptionnel. 
{\em Topology} {\bf 26} (1987), 93-105.\\

\bibitem{GS} {\sc Ghys, \'E. \& Sergiescu, V.} Sur un groupe remarquable de
diff\'eomorphismes 
du cercle. {\em Comment. Math. Helv.} {\bf 62} (1987), 185-239.\\

\bibitem{GT} {\sc Ghys, \'E. \& Tsuboi, T.} Diff\'erentiabilit\'e des conjugaisons
entre 
syst\`emes dynamiques de dimension 1. {\em Annales de l'Institut Fourier (Grenoble)} 
{\bf 38} (1988), 215-244.\\

\bibitem{imca} {\sc Navas, A.} {\em Grupos de difeomorfismos del c\'{\i}rculo.} 
Monograf\'{\i}as del IMCA, Lima, Per\'u (2006).\\

\bibitem{quasim} {\sc Navas, A.} On uniformly quasisymmetric groups of 
circle diffeomorphisms. {\em An. Acad. Sci. Fenn. Math.} {\bf 31} 
(2006), 437-462.\\

\bibitem{brasil} {\sc Oliveira, F. \& da Rocha, L.} Minimal non ergodic
$C^1$-diffeomorphisms 
of the circle. {\em Erg. Theory. and Dyn. Systems} {\bf 21} (2001), 1843-1854.\\

\bibitem{quas} {\sc Quas, A.} Non-ergodicity for $C^1$ expanding maps and
$g$-measures. 
{\em Erg. Theory and Dyn. Systems} {\bf 16} (1996), 531-543.\\

\end{thebibliography}
\end{document}